\documentclass[12pt,leqno]{amsart}
\usepackage{amssymb}
\usepackage{amsmath,amssymb,color}
\begin{document}
\baselineskip=18pt

\newtheorem{theorem}{Theorem}

\newtheorem{defi}{Definition}

\newtheorem{corollary}[theorem]{Corollary}

\newtheorem{proposition}{Proposition}

\newtheorem{lemma}{Lemma}

\newtheorem{remark}{Remark}


%
\newtheorem{thm}{Theorem}[section]
\newtheorem{lem}[thm]{Lemma}
\newtheorem{cor}[thm]{Corollary}
\newtheorem{prop}[thm]{Proposition}
\newtheorem{defn}[thm]{Definiction}
\newtheorem{rmk}{Remark}
\newtheorem{exa}{Example}
\numberwithin{equation}{section}
\newcommand{\DDDD}{D_{\rho,\delta}}
\newcommand{\weight}{e^{2s\varphi}}
\newcommand{\tilweight}{e^{2s\widetilde{\va}}}
\newcommand{\ep}{\varepsilon}
\newcommand{\la}{\lambda}
\newcommand{\va}{\varphi}
\newcommand{\ppp}{\partial}
\newcommand{\izt}{\int^t_0}
\newcommand{\ddd}{\mathcal{D}}
\newcommand{\www}{\widetilde}
\newcommand{\pppg}{\partial_{t,\gamma}^{\alpha}}
\newcommand{\pppa}{\partial_t^{\alpha}}
\newcommand{\VVVV}{Cs^3e^{Cs}\mathcal{D}(u)^2 + Cs^3e^{2s\sigma_1}M^2}
\newcommand{\ddda}{d_t^{\alpha}}
\newcommand{\DDDa}{D_t^{\alpha}}
\newcommand{\lla}{L^{\frac{3}{4}}}
\newcommand{\llb}{L^{-\frac{1}{4}}}
\newcommand{\HHMIN}{H^{-1}(\Omega)}
\newcommand{\HHONE}{H^1_0(\Omega)}
\newcommand{\sumij}{\sum_{i,j=1}^n}
\newcommand{\Halp}{H_{\alpha}(0,T)}
\newcommand{\Hal}{H^{\alpha}(0,T)}
\newcommand{\lllb}{L^{-\frac{3}{4}}}
\newcommand{\pdif}[2]{\frac{\partial #1}{\partial #2}}
\newcommand{\ppdif}[2]{\frac{\partial^2 #1}{{\partial #2}^2}}
\newcommand{\uen}{u_N^{\ep}}
\newcommand{\fen}{F_N^{\ep}}
\newcommand{\penk}{p^{\ep}_{N,k}}
\newcommand{\penl}{p^{\ep}_{N,\ell}}
\newcommand{\R}{\mathbb{R}}
\newcommand{\Q}{\mathbb{Q}}
\newcommand{\Z}{\mathbb{Z}}
\newcommand{\C}{\mathbb{C}}
\newcommand{\N}{\mathbb{N}}
\newcommand{\uu}{\mathbf{u}}
\renewcommand{\v}{\mathbf{v}}
\newcommand{\y}{\mathbf{y}}
\newcommand{\RR}{\mathbf{R}}
\newcommand{\Y}{\mathbf{Y}}
\newcommand{\w}{\mathbf{w}}
\newcommand{\z}{\mathbf{z}}
\newcommand{\G}{\mathbf{G}}
\newcommand{\ooo}{\overline}
\newcommand{\OOO}{\Omega}
\newcommand{\CCO}{{_{0}C^1[0,T]}}
\newcommand{\dddaO}{{_{0}\ddda}}
\newcommand{\WWW}{{_{0}W^{1,1}(0,T)}}
\newcommand{\CCOO}{{^{0}C^1[0,T]}}
\newcommand{\HH}{_{0}{H^{\alpha}}(0,T)}
\newcommand{\HHO}{_{0}H^1(0,T)}
\newcommand{\hhalf}{\frac{1}{2}}
\newcommand{\DDD}{\mathcal{D}}
\newcommand{\RRR}{\mathcal{R}}
\newcommand{\RRRR}{\longrightarrow}
\newcommand{\lbra}{_{H^{-1}(\Omega)}{\langle}}
\newcommand{\wdel}{\www{\delta}}
\newcommand{\weps}{\www{\varepsilon}}
\newcommand{\rbra}{\rangle_{H^1_0(\Omega)}}
\newcommand{\lbran}{_{(H^{-1}(\Omega))^N}{\langle}}
\newcommand{\rbran}{\rangle_{(H^1_0(\Omega))^N}}
\newcommand{\llll}{L^{\infty}(\Omega\times (0,t_1))}
\newcommand{\sumk}{\sum_{k=1}^{\infty}}
\renewcommand{\baselinestretch}{1.5}
\renewcommand{\div}{\mathrm{div}\,}  
\newcommand{\grad}{\mathrm{grad}\,}  
\newcommand{\rot}{\mathrm{rot}\,}  
\newcommand{\xxx}{x', x_n}
\newcommand{\xxz}{x', 0}
\newcommand{\urur}{\left( \frac{1}{r}, \theta\right)}

\allowdisplaybreaks
\renewcommand{\theenumi}{\arabic{enumi}}
\renewcommand{\labelenumi}{(\theenumi)}
\renewcommand{\theenumii}{\alph{enumii}}
\renewcommand{\labelenumii}{(\theenumii)}
\def\thefootnote{{}}

\title
[]
{Sharp uniqueness and stability of solution for  an inverse source problem
for the Schr\"odinger equation
}

\author{
$^1$ O.Y. Imanuvilov and 
$^{2}$ M.~Yamamoto }

\thanks{
$^1$  
Department of Mathematics, Colorado State University, 101 Weber Building, 
Fort Collins CO 80523-1874, USA e-mail:{\tt oleg@math.colostate.edu}
\\
$^2$ Graduate School of Mathematical Sciences, The University
of Tokyo, Komaba, Meguro, Tokyo 153-8914, Japan 
e-mail: {\tt myama@ms.u-tokyo.ac.jp}
}

\date{}
\maketitle

\baselineskip 18pt

\begin{abstract}
 The manuscript is concerned with  uniqueness and stability for inverse source 
problem of determining spatially varying factor $f(x)$ of a source term 
given by $R(t)f(x)$ with suitable given $R(t)$ in the right hand side of the  Schr\"odinger equation with time independent coefficients.
In order to establish these results 
we provide a simple proof of a logarithmic conditional stability  of  the Cauchy problem  
for the  Schr\"odinger equation with time-independent coefficients and the zero Dirichlet boundary conditions  on the whole boundary and a proof of uniqueness of solution to the Cauchy problem for the Schr\"odinger equation with data on an arbitrary small part of a lateral boundary.
We do not assume any geometrical constraints on subboundary 
and a time interval is arbitrary.  
The key is an integral transform, with a kernel solving a null controllability problem for the 1-D Schr\"odinger,   which changes a solution 
of the Schr\"odinger equation to a solution of an elliptic equation.

\end{abstract}

\date{}

\baselineskip 18pt


\section{Introduction and main results}

Let $\Omega$ be a bounded domain in $\Bbb R^n$ with $C^2$-boundary.
We set $x= (x_1, ..., x_n) \in \R^n$, $\ppp_tu = \frac{\ppp u}{\ppp t}$,
$\ppp_{x_j}u = \frac{\ppp u}{\ppp x_j}$ for $1\le j \le n$, and
$\nabla = (\ppp_{x_1}, ..., \ppp_{x_n})$.
Moreover, $\nu = \nu(x)$ denotes the unit outward normal vector to $\ppp\OOO$
at $x$ and $\ppp_{\nu}u = \nabla u \cdot \nu$ on $\ppp\OOO$.
We set $Q:= (0,T) \times \OOO$ and $i:= \sqrt{-1}$, $D=(D_0,D'), D'=(D_1,\dots, D_n), D_0=\frac 1i\partial_t, D_j=\frac 1i\partial_{x_j}, j\in \{1,\dots, n\}.$

We consider the  Schr\"odinger equation:
\begin{equation}\label{(1.1)}
P(x,D)u:= i\partial_tu-\sum_{j,k=1}^n \partial_{x_j}(a_{jk}(x)\partial_{x_k} u(t,x))
+ \sum_{j=1}^n b_j(x)\partial_{x_j} u(t,x) 
+ c(x)u(t,x),\quad (t,x) \in Q,                  
\end{equation}

Throughout this article, we make the following assumptions: 
the coefficients $a_{jk}$ 
are real-valued  functions for all $1\le j,k \le n$ and 
\begin{equation}\label{SKP}
a_{jk}=a_{kj} \in C^{1}(\ooo{\OOO}), \quad b_j\in  W^{1}_\infty(\Omega), \quad c \in L^{\infty}(\OOO) 
\quad \mbox{for all $j,k \in \{1,\dots,n\}$}    
\end{equation}
and there exists a constant $\beta>0$ such that
\begin{equation}\label{SKP1}
\sum_{j,k=1}^n a_{jk}(x)\xi_j\xi_k \ge \beta\sum_{j=1}^n \xi_j^2 \quad
\mbox{for all $x \in \OOO$, $\xi_1, ..., \xi_n \in \R$}.     
\end{equation}

We set 
$$
A(x,D')v:= -\sum_{j,k=1}^n \partial_{x_j}(a_{jk}(x)\partial_{x_k} v(x))
+ \sum_{j=1}^n b_j(x)\partial_{x_j} v(x) + c(x)v(x)
$$
for $x \in \OOO$.

Let $\Gamma \subset \ppp\OOO$ be a suboundary and let $T>0$.

We define
$$
H^{1,2}(Q):= \{ u \in L^2(Q);\,
\ppp_tu,\, \ppp_x^{\alpha}u \in L^2(Q) \,\, \mbox{if $\vert \alpha
\vert \le 2$} \},
$$
where $\alpha:= (\alpha_1, ..., \alpha_n) \in (\N \cup \{0\})^n$,
$\vert \alpha\vert := \alpha_1 + \cdots + \alpha_n$ and 
$\ppp_x^{\alpha}:= \ppp_{x_1}^{\alpha_1} \cdots \ppp_{x_n}^{\alpha_n}$.

In this manuscript we concern with the following inverse problem:
{\it 

Let $u \in C([0,T];L^2(\OOO))$ and $f \in L^2(\OOO)$ satisfy
\begin{equation}\label{(1.7)}
P(x,D)u= i\partial_tu + A(x,D')u = R(t)f(x)\quad \mbox{in}\quad Q,\quad u(0,\cdot)=0. 
\end{equation}
Here $R$ is a given real-valued  function. Determine the  source term $f$ based on the following data : $\partial_\nu u\vert_\Gamma.$}

We have 
\\
\begin{theorem}\label{3} Let pair  $(u,f)$ satisfies (\ref{(1.7)}) and $R\in C^1[0,T].$\newline
{\it
(i) 
We assume that $u,\partial_t u\in H^{1,2}(Q)$ and
\begin{equation}\label{(1.8)}
R(0) \ne 0.                
\end{equation}
If $u = \ppp_{\nu}u = 0$ on $(0,T) \times \Gamma$ with 
arbitrary fixed subdomain $\Gamma$ of $\partial\Omega$, then 
$f = 0$ in $\OOO$ and $u=0$ in $Q$.
\\
(ii) We assume that $u\in C([0,T];L^2(\Omega)), f\in L^2(\Omega)$ and
\begin{equation}\label{soplo}
\mbox{$R$ is not identically zero on $[0,T]$.}   
\end{equation}
Let $u =0$ on $(0,T) \times \ppp\OOO$ and
$\ppp_{\nu}u = 0$ on $(0,T) \times \Gamma$ with 
arbitrary fixed subdomain $\Gamma$ of $\partial\Omega$, then 
$f = 0$ in $\OOO$ and $u=0$ in $Q$.
}
\end{theorem}

{\bf Remark.} {\it  A function $u \in C([0,T];L^2(\Omega))$  is a 
solution (\ref{(1.7)}) such that  $u\vert_{(0,T)\times \partial\Omega} 
= \ppp_{\nu}u\vert_{(0,T) \times \Gamma} = 0$ 
if 
$$
(u,P^*(x,D)w)_{L^2(Q)}=(Rf, w)_{L^2(Q)} \quad \forall w\in H^{1,2}(Q), w=0\quad \mbox{on}\,\, (0,T)\times  \partial\Omega\setminus\Gamma, \, w(T,\cdot)=0.
$$}

We do not know whether $f=0$ in $\OOO$ and $u=0$ in $Q$ follow from 
$u = \ppp_{\nu}u = 0$ on $(0,T) \times \Gamma$ under the weaker 
assumption (\ref{soplo}). 
Moreover, in terms of Theorem \ref{3}, we can prove
\\
\begin{theorem}\label{4}
{\it Let pair  $(u,f)\in H^{1,2}(Q)\times L^2(\Omega)$ satisfies (\ref{(1.7)}), $R\in C^2[0,T]$. Moreover $\partial^2_tu, \partial_t u\in H^{1,2}(Q)$ and $u\vert_{(0,T)\times \partial\Omega}=0.$
We assume (\ref{(1.8)}) and 
$$
\Vert \ppp_tu\Vert_{H^{1,2}(Q)}
+ \Vert u\Vert_{H^{1,2}(Q)} \le M,
$$
where $M>0$ is an arbitrarily chosen constant.
Then there exists a constant $C>0$ such that
\begin{equation}\label{cocos}
\Vert f\Vert_{L^2(\Omega)}
\le \frac{C}{\vert \log \Vert \partial_\nu u\Vert
_{H^2(0,T;L^2(\Gamma))}\vert}.
\end{equation}
}
\end{theorem}
The inverse problem of  recovery of a source term in the right hand side of the Schr\"odinger equation  was studied in 
Baudoin and Mercado \cite{BM}, Baudouin and Puel \cite{B-P}, \cite{BP1},
Mercado, Osses and Rosier \cite{MOR}, Yuan and Yamamoto \cite{YM}.

The authors, in the above mentioned papers, studied right hand side of more general form, namely $R(t,x)f(x)$ and proved uniqueness of determination of the right hand side  and  Lipschitz stability of the determination. The machinery of these papers  is based on the Carleman estimates with boundary for the Schr\"odinger  equation.
Such an approach is based on pioneering paper of A. Bukhgeim and M. Klibanov \cite{BK}. Unfortunately in the case of Schr\"odinger equation a Carleman estimate  with boundary require an existence of an appropriate pseudoconvex function. This condition is  restrictive. For example, in case of the principal part of $A$ be the Laplace operator, it require for observed part of the boundary $\Gamma$ to be sufficiently large. In the case of the general elliptic operator $A$ existence of pseudoconvex function implies the non-trapping condition for the Hamiltonian flow with the Hamiltonian $H(x,\xi)= \sum_{j,k=1}^n a_{jk}(x)\xi_j\xi_k.$ 

Proof of the theorem \ref{3} and theorem \ref{4}  are based on the following uniqueness and conditional stability results 
 for the Schr\"odinger equation:

Assuming the zero Dirichlet boundary condition on 
the whole lateral boundary $(0,T) \times \ppp\OOO$, we can prove 
a conditional stability estimate of logarithmic rate.
\begin{theorem}\label{2}
{\it 
Let a suboundary $\Gamma \subset \ppp\OOO$ and a constant $M>0$ 
be arbitrarily chosen.
Let $u,\partial_t u \in H^{1,2}(Q)$ satisfy 
\begin{equation}\label{(1.4)}
P(x,D) u: =i\partial_tu + A(x,D')u=0\quad \mbox{in}\quad Q
                                             \end{equation}

and
\begin{equation}\label{(1.5)}
u=0 \quad \mbox{on $(0,T)\times \partial\Omega$}.          
\end{equation}
If 
\begin{equation}\label{(1.6)}
\Vert u\Vert_{H^{1,2}((0,T)\times \Omega)}+\Vert\partial_t u\Vert_{H^{1,2}((0,T)\times \Omega)} \le M,       
\end{equation}
then there exists a constant $C>0$  independent of $u$ such that
\begin{equation}\label{lopukh5}
\Vert u\Vert_{ C([0,T];L^2(\Omega))}
\le \frac{C}{\vert \log \Vert \partial_\nu u\Vert_{H^1((0,T)\times \Gamma)}
\vert}.
\end{equation}}
\end{theorem}

We emphasize that in addition to (\ref{SKP}) and (\ref{SKP1}),
we do not assume any extra conditions for the principal 
coefficients $a_{jk}$, no geometrical constraints on $\Gamma$ and no
limitation of length of time interval of observation $T.$

We remark that if $a_{jk}$ satisfies some conditions called the 
pseudo-convexity (e.g., \cite{B-P}, \cite{BP1},  
Yuan and Yamamoto \cite{YM}) and $\Gamma$ is some large portion of 
$\ppp\OOO$, then one can prove the Lipschitz stability rather than 
the logarithmic stability.  In a particular case where 
$a_{jj} = 1$ and $a_{jk} = 0$ for $j\ne k$ and $j,k=1,..., n$,
if $\Gamma \supset \{x\in \ppp\OOO;\, (x-x_0)\cdot\nu(x) \ge 0\}$
with some $x_0\in \R^n$ (e.g., \cite{B-P}, \cite{BP1}), then
$$
\Vert u(0,\cdot)\Vert_{L^2(\OOO)} \le 
C\Vert \ppp_{\nu}u\Vert_{L^2((0,T)\times \Gamma)}.
$$
Here we emphasize that this is unconditional stability without 
any a priori condition such as (\ref{(1.6)}).

The conclusion holds under the a priori boundedness condition (\ref{(1.6)}) and so 
is conditional stability.

Finally we  formulate the uniqueness result  for the Cauchy problem for the Schr\"odinger equation with solution from 
$C([0,T];L^2(\Omega)).$
\\
\begin{theorem}\label{1}
{\it 
Let $\Gamma \subset \ppp\OOO$ be an arbitrarily chosen subboundary 
and $T>0$ be arbitrary.  
We assume that $u \in C([0,T];L^2(\Omega))$ satisfies (\ref{(1.1)}).
Then $u = \ppp_{\nu}u = 0$ on 
$(0,T) \times \Gamma$ yields $u=0$ in $(0,T) \times \OOO$.
}
\end{theorem}

{\bf Remark.} {\it  Function $u \in C([0,T];L^2(\Omega))$  is solution (\ref{(1.1)}) such that  $u = \ppp_{\nu}u = 0$ on 
$(0,T) \times \Gamma$ if 
$$
(u,P^*(x,D)w)_{L^2(Q)}=0 \quad \forall w\in H^{1,2}(Q),$$
$$ w=\partial_\nu w=0\quad \mbox{on}\,\, (0,T)\times  \partial\Omega\setminus\Gamma,\, w(T,\cdot)=w(0,\cdot)=0.
$$}

The unique continuation is a fundamental property of partial differential 
equations.  In order to conclude that $u=0$ in $Q$, the conditions on 
$\Gamma\subset \ppp\OOO$ and $T>0$ depend on the types of partial 
differential equations under consideration.
For the parabolic equation $\ppp_tu + Au = 0$, for any subboundary $\Gamma$ and
any $T>0$, without any assumptions on $a_{jk}$ except for (\ref{SKP}), (\ref{SKP1}) we can 
prove the unique continuation (e.g., Mizohata \cite{Mi}, Saut and
Scheurer \cite{SS}).  Here we do not intend to create any comprehensive 
references.  Moreover, the same sharp unique continuation holds also in the 
case where the coefficients depend on $t$ and satisfy (\ref{SKP}) and (\ref{SKP1}).
Here and henceforth 
we do not exploit the minimum regularity of the coefficients.   

For the hyperbolic equation $\ppp_t^2u + Au = 0$, in the case where 
$a_{jk}$ do not depend on $t$ or analytic in $t$, a sharp unique continuation 
of global type similar to the parabolic equation holds provided that 
$T>0$ is sufficiently large (e.g., Robbiano \cite{Luc}, Tataru \cite{Tat}).

Compared with the parabolic and hyperbolic equations, it seems that 
the researches on the unique continuation for Schr\"odinger equations,
are not completed.  Although we may be suggested to be able to prove
the sharp unique continuation similar to the parabolic equation by 
comprehensive researches on Carleman estimates, we can not find 
a short proof in the existing publications, and one of  main purposes of this 
article is to provide such a self-contained and simple proof.

The logarithmic conditional stability results  for the unique continuation of the hyperbolic equation were established in \cite{BO1} and \cite{BO2}.

The article is composed of six sections and appendix.  In Sections 2 and 3, we prove
Theorem \ref{2} statement I and Theorem  \ref{1} respectively.  Section 4 is devoted to the proof Theorem \ref{2} statement II.   In Section 5, prove Theorem \ref{3} and in section 6 we give a proof of Theorem \ref{4}. In appendix we solve some controllability problem for 1-D Schr\"odinger equation in order to construct a kernel of integral transform.
\section{Proof of Theorem \ref{3}, statement  I}

We define the operator by 
\begin{equation}\label{(4.1)}
(Kv)(t):= R(0) v(t) + \int^t_0 R'(t-s)v(s) ds, \quad 0<t<T.
                                                                 \end{equation}
Let $u$ satisfy (\ref{(1.7)}) and $u(0,\cdot) = 0$ in $\OOO$ and
$u = \vert \nabla u\vert = 0$ on $(0,T)\times \Gamma$.  
Consider the  equation with respect to $z(t,x)$:
\begin{equation}\label{(4.2)}
\ppp_tu(t,x) = (Kz)(t,x), \quad 0<t<T, \, x\in \OOO.
                                                        \end{equation}

Since $R(0) \ne 0$, the operator $K$ is the Volterra operator of the second 
kind, and we see that $K^{-1}: H^1(0,T) \longrightarrow H^1(0,T)$ exists and 
is bounded.  Therefore, $z(\cdot,x) \in H^1(0,T)$ is well defined for each 
$x \in \OOO$, and
\begin{equation}\label{(4.3)}
\ppp_tu(t,x) = R(0)z(t,x) + \int^t_0 R'(t-s)z(s,x) ds, \quad 0<t<T,\, x\in 
\OOO                                                       
\end{equation}
by $\ppp_tu \in H^{1,2}(Q)$.

Since $u(0,\cdot) = 0$ in $\OOO$ and 
$$
R(0)z(t,x) + \int^t_0 R'(t-s)z(s,x) ds
= \ppp_t\left( \int^t_0 R(t-s)z(s,x) ds\right), 
$$
we obtain
$$u(t,x) = \int^t_0 R(t-s)z(s,x) ds, \quad 0<t<T,\, x\in \OOO,
$$
that is,
\begin{equation}\label{(4.4)}
u(t,x) = \int^t_0 R(s)z(t-s,x) ds, \quad 0<t<T,\, x\in \OOO.
                                                                \end{equation}

We will prove that $z \in H^{1,2}(Q)$ satisfies 
\begin{equation}\label{(4.5)}
i\ppp_tz(t,x) + A(x,D')z(t,x)=0 \quad \mbox{in $(0, T) \times \OOO$},
                                                                \end{equation}
where $t_*\in (0,T)$ is some constant, and 
\begin{equation}\label{(4.6)}
z = \vert \nabla z\vert = 0 \quad \mbox{on $(0,T) \times \Gamma$}.
                                                           \end{equation}

We can readily verify (\ref{(4.6)}), because $\ppp_tu = \vert \nabla\ppp_tu\vert 
= 0$ on $(0,T) \times \Gamma$ implies 
$Kz = \vert K(\nabla z)\vert = 0$ on $(0,T) \times \Gamma$, so that  
the injectivity of $K$ directly yields (\ref{(4.6)}).

Using $\ppp_tu \in H^{1,2}(Q) \subset C([0,T];L^2(\OOO))$ by (\ref{(4.3)}), we have 
$$
\ppp_tu(0,x) = R(0)z(0,x), \quad x\in \OOO.
$$
On the other hand, substituting $t=0$ in (\ref{(1.7)})  we obtain
$$
\ppp_tu(0,x) = -iR(0)f(x), \quad x\in \OOO.
$$
Hence $R(0)z(0,x) = -iR(0)f(x)$ for $x\in \OOO$.  By $R(0) \ne 0$, we 
reach 
\begin{equation}\label{(4.7)}
z(0,x) = -if(x), \quad x\in \OOO.                  \end{equation}

Now we will prove (\ref{(4.5)}).  In terms of (\ref{(4.4)}) and (\ref{(4.7)}), we have
\begin{align*}
& \ppp_tu(t,x) = R(t)z(0,x) + \int^t_0 R(s)\ppp_tz(t-s,x) ds\\
=&-i R(t)f(x) + \int^t_0 R(s)\ppp_tz(t-s,x) ds
\end{align*}
and
$$
A(x,D')u(t,x) = \int^t_0 R(s)A(x,D')z(t-s,x) ds, \quad 0<t<T,\, x\in \OOO.
$$
Consequently (\ref{(1.7)}) implies
\begin{align*}
& -iR(t)f(x) = \ppp_tu - iA(x,D')u(t,x)\\
=& -iR(t)f(x) + \int^t_0 R(s)(\ppp_tz - iA(x,D')z)(t-s,x) ds,
\end{align*}
that is,
$$
\int^t_0 R(s)(\ppp_tz -i A(x,D')z)(t-s,x) ds = 0, \quad 0<t<T,\, x\in \OOO.
$$
Hence, setting $Z(s):= \Vert (\ppp_sz - iA(x,D')z)(s,\cdot)\Vert_{L^2(\OOO)}$ 
for $0<t<T$, we reach 
$$
\int^t_0 R(s) Z(t-s) ds = 0, \quad 0<t<T.
$$
By the Titchmarsh convolution  theorem (e.g., Titchmarsh \cite{Tit}), 
there exists $t_* \in [0,T]$ such that 
$$
R(s) = 0 \quad \mbox{for $0<s<T-t_*$}, \quad
Z(s) = 0 \quad \mbox{for $0<s<t_*$}.
$$
Since $R(0)\ne  0$ in $[0,T]$, we see that $t^*=T.$  Thus the verification of (\ref{(4.5)}) is complete.

From (\ref{(4.5)}) and (\ref{(4.6)}) applying  the uniqueness theorem \ref{1} we obtain $z(0,\cdot)=0.$ Then equality (\ref{(4.7)}) implies $f=0.$
$\blacksquare$

\section{Proof of Theorem \ref{1}. }\label{ZZZQ!88}


In order to prove the second part of theorem \ref{3}  we need  the uniqueness result  for the Cauchy problem for the Schr\"odinger equation with data on an arbitrary small part of the boundary.

{\bf Proof.} Let $\widetilde \Omega$ be a bounded domain in 
$\Bbb R^n$ with smooth boundary  such that $\Omega\subset \widetilde \Omega$,
$\partial\Omega\setminus\Gamma\subset \partial\widetilde \Omega$ and 
$\omega=\widetilde \Omega\setminus\overline\Omega$ is an open set.  We extend a function $u$ by zero on $(0,T)\times \tilde \Omega\setminus \Omega$. Function $u$ satisfies 
\begin{equation}\label{Simas} P(x,D)u=0\quad \mbox{in}\,\, (0,T)\times \tilde\Omega,\quad u\vert_{(0,T)\times \omega}=0.\end{equation}
Let $\mathcal K(t,\tau)\in C^\infty([-1,1]\times[0,T])$ be solution to the boundary controllability problem for the Shr{\"o}dinger equation:
\begin{eqnarray}i\partial_{\tau}\mathcal K-\partial^2_{t}\mathcal K =0\quad  \mbox{in} \quad  t\in (-1,1)\,\,\tau\in (0,T),\nonumber\\
\quad \mathcal K(-1,\tau)=\psi (\tau)\,\,\mbox{on}\,\, (0,T),
\,\quad \mathcal K(\cdot,0)=\mathcal K(\cdot, T)=0.\end{eqnarray} Here $\psi\in C^\infty_0(0,T)$ is fixed complex valued function. Control is located on the part of the boundary $t=1.$ For existence of such a solution see e.g. \cite{MRR}.
Consider the function $w(x)=\int_{0}^T\mathcal K(t,\tau)u(\tau,x)d\tau.$ This function satisfies the equation
\begin{equation}\label{simas}
L(x,D) w=-\partial_{t}^2w+A(x,D')w=0 \quad \mbox{in}\quad  (-1,1)\times \tilde \Omega,\quad w=0\quad \mbox{on}\quad (-1,1)\times \omega.
\end{equation}
The second inequality in (\ref{simas}) follows from (\ref{Simas}) immediately.
In order to prove the first equality we observe that for any $p\in H^2_0((-1,1)\times \tilde \Omega)$
$$
(w,L^*(x,D)p)_{L^2((-1,1)\times \tilde \Omega)}=(w,-\partial_{t}^2p+A^*(x,D')p)_{L^2((-1,1)\times \tilde \Omega)}=
$$
$$(\int_{0}^T\mathcal K(t,\tau)u(\tau,x) d\tau,(-\partial_{t}^2p+A^*(x,D')p)_{L^2((-1,1)\times \tilde \Omega)}=
$$
$$
(u,-\int_{-1}^1\mathcal K(t,\tau)\partial_{t}^2pdt+A^*(x,D')\int_{-1}^1\mathcal K(t,\tau)p dt)_{L^2((0,T)\times \tilde \Omega)}=
$$
$$
(u,-i\partial_{\tau}\int_{-1}^1\mathcal K(t,\tau)pdt+A^*(x,D')\int_{-1}^1\mathcal K(t,\tau)p dt)_{L^2((0,T)\times \tilde \Omega)}=(u,P^*(x,D)\tilde p)_{L^2(Q)}=0
$$
where 
$\tilde p(\tau,x)=\int_{-1}^{1}\mathcal K(t,\tau)p(t,x)dt.$

The equations (\ref{simas}) are established. By assumptions  (\ref{SKP}), (\ref{SKP1}) the operator $L$ is the second order elliptic operator with $C^1$ coefficients in the principal part and $L^\infty$ coefficients in first and zero order terms. Then theorem  8.9.1 from \cite{Hor}  which establish uniqueness of the Cauchy problem for  this operator can be applied. This implies that function $w$ is identically equal to zero on the cylinder $(-1,1)\times \Omega.$ In particular
$$
w(-1,x)=\int_{0}^T\mathcal K(-1,\tau)u(\tau,x)d\tau=\int_{0}^T \psi(\tau)u(\tau,x)d\tau=0\quad \mbox{on}\quad \Omega.
$$
Since the space $C^\infty_0(0,T)$ is dense in $L^2(0,T)$ we have $u\equiv 0.$ 
Proof of theorem is complete. $\blacksquare$

\section{Proof of Theorem \ref{3}, statement II}

{\bf Proof.}
Without loss of generality we may assume that \begin{equation}\label{soplo10}0\in \mbox{supp}\,R.\end{equation} Indeed, since $\mbox{supp}\,R$ is not empty there exist a point $\tilde t\in \mbox{supp}\,R. $ Denote by $t_0$ infimum over such points. Function $u$ is zero on $(0,t_0)\times \Omega.$   We replace function $u$ and $R$  in (\ref{(1.7)}) by $y(t+t_0,\cdot)$ and $R(t+t_0)$ respectively to obtain (\ref{soplo10}).

Next we reduce the inverse  problem to the case when $R(0)=0.$ Function $w(x)=\int_0^{t}u(\tau,x)d\tau$ satisfies :
$$
P(x,D) w=\int_0^{t}R(\tau)d\tau f(x)\quad \mbox{in}\quad Q,\quad  w(0,\cdot)=0, \quad w\vert_{(0,T)\times \partial\Omega}=\partial_\nu w\vert_{(0,T)\times \Gamma}=0.
$$ Obviously function $\tilde R(t)=\int_0^t R(\tau)d\tau$ satisfies 
$$
\tilde R(0)=0.$$ Moreover, since $\frac{d\tilde R}{dt}=R$ by (\ref{soplo}), function $\tilde R$ is not identically equal zero on $(0,T).$ 
Hence, besides (\ref{soplo})
we may additionally assume 
\begin{equation}\label{fuk}R(0)=0.\end{equation} 
Let $\widetilde \Omega$ be a bounded domain in 
$\Bbb R^n$ with smooth boundary  such that $\Omega\subset \widetilde \Omega$,
$\partial\Omega\setminus\Gamma\subset \partial\widetilde \Omega$ and 
$\omega=\widetilde \Omega\setminus\overline\Omega$ is an open set.  We extend  functions $u$  and $f$ by zero on $(0,T)\times \tilde \Omega\setminus \Omega$. 
By our assumptions we have 
$$ P(x,D)u=Rf\quad \mbox{in}\quad (0,T)\times \tilde \Omega, \quad u\vert_{(0,T)\times \partial\tilde \Omega}=0,\,\,
u=0\quad \mbox{on}\quad   (0,T)\times \omega.$$
Next we consider function 
\begin{equation}\label{popkorn}
y(t,x)=\int_0^{t}R(t-\tau)v(\tau,x)d\tau, \quad x\in \tilde \Omega,
\end{equation}
where function $v\in C([0,T];L^2(\tilde\Omega))$ solves the initial value problem
$$
P(x,D)v=0\quad \mbox{in}\quad (0,T)\times \tilde \Omega,\quad  v(0,\cdot)=-if, \quad v\vert_{(0,T)\times \partial\tilde \Omega}=0.
$$
By (\ref{popkorn})
function $y(t,x)$ satisfies the initial and boundary conditions 
\begin{equation}\label{lop}\quad y\vert_{(0,T)\times \partial\tilde \Omega}=0,\quad  y(0,\cdot)=0.
\end{equation}
The short computations imply
\begin{equation}\label{lopp}
P(x,D)y=Rf\quad \mbox{in}\quad (0,T)\times \tilde \Omega.
\end{equation}
Indeed, by (\ref{fuk}) for any  $p\in H^{1,2}((0,T)\times \tilde \Omega), p(T,\cdot)=0, p\vert_{(0,T)\times \partial\tilde \Omega}=0$ we have  $$ (y,P^*p)_{L^2((0,T)\times \tilde \Omega)}=(y, i\partial_t p +A^*p)_{L^2((0,T)\times \tilde \Omega)}=
(\int_0^{t}R(t-\tau)v(\tau,x)d\tau, i\partial_t p +A^*p)_{L^2((0,T)\times \tilde \Omega)}=
$$
$$(\int_0^{t}R(t-\tau)v(\tau,x)d\tau, A^*p)_{L^2((0,T)\times \tilde \Omega)}-(\int_0^{t}\partial_tR(t-\tau)v(\tau,x)d\tau, ip )_{L^2((0,T)\times \tilde \Omega)}=
$$
$$(\int_0^{t}R(t-\tau)v(\tau,x)d\tau, A^*p)_{L^2((0,T)\times \tilde \Omega)}+(\int_0^{t}\partial_\tau R(t-\tau)v(\tau,x)d\tau, ip )_{L^2((0,T)\times \tilde \Omega)}=
$$
$$\int_0^{T}(v, P^* R(t-\cdot) p(t,\cdot))_{L^2((0,T)\times \tilde \Omega)}dt=(f,p)_{L^2((0,T)\times \tilde \Omega)}.
$$
Hence function $y$ solves the initial value problem (\ref{lopp}), (\ref{lop}). By the uniqueness of solution to this initial value problem $u=y$ and 
$$ y\vert_{(0,T)\times \omega}=0.$$

This equality implies 
$$
\int_0^{t}R(t-\tau) v(\tau,x)d\tau=0\quad\mbox{on}\quad  (0,T)\times \omega.
$$
Hence 
$$
\int_0^{t}R(t-\tau) \Vert v(\tau,\cdot )\Vert_{L^2(\omega)}d\tau=0.
$$
 
Applying the Titchmarsh convolution  theorem  we have there exist $\lambda\ge 0$ and $r\ge 0$ such that
\begin{equation}\label{oplot}
\lambda+r\ge T
\end{equation}
and function $R$ is equal to zero on $(0,\lambda)$ and   function $\Vert v(\tau,\cdot )\Vert_{L^2(\omega)}$ is equal to zero on $(0,r). $ By (\ref{soplo}) and (\ref{soplo10}) parameter $\lambda=0.$  Hence  by (\ref{oplot}) we have $r=T$ and 
$$
v=0\quad \mbox{on}\quad   (0,T)\times \omega.$$
Then theorem \ref{1} implies 
$$
v(t,x)=0 \quad \mbox{on}\quad (0,T)\times \Omega.
$$
Hence, by (\ref{popkorn}) $$ y(t,x)=0 \quad \mbox{on}\quad (0,T)\times \Omega.$$ This implies 
$$
Rf=0\quad \mbox{on}\quad (0,T)\times \Omega.$$ Since function $R$ is not identically equal zero on $(0,T)$ we have $f\equiv 0.$ $\blacksquare$

\section{Proof of Theorem \ref{2}.}

{\bf Proof.}  Without loss of generality we may assume that 
$$
\Vert \partial_\nu u\Vert_{L^2((0,T)\times\Gamma)}\le 1\quad \mbox{and}\quad \Vert \partial_\nu\partial_t u\Vert_{L^2((0,T)\times\Gamma)}\le 1.
$$
Let $y\in\Omega$ and  $\mathcal K(t,\tau)$ be solution to the boundary controllability problem for the one dimensional  Shr{\"o}dinger equation:
\begin{eqnarray}\label{HUKK}i\partial_{\tau}\mathcal K-\partial^2_{t}\mathcal K =0\quad  \mbox{in} \quad  t\in (0,2), \,\tau\in(0,T),\nonumber\\
\quad \mathcal K(0,\tau)=\psi(\tau)\,\,\mbox{on}\,\, (0,T),
\quad \mathcal K(\cdot,0)=\mathcal K(\cdot,T)=0.\end{eqnarray} 

Control is located on the part of the boundary $t=2:$
$$
\mathcal K(2,\tau)= v(\tau).
$$ A given  function $\psi$ is smooth function but in general does not satisfy to the compatibility conditions at  time moment $\tau=0$ or $\tau=T.$ Solution to the controllability problem (\ref{HUKK}) is not unique and in general is not smooth, since  function  $\psi$ in general does not satisfy the compatibility conditions. On the other hand one can chose solution to the controllability problem in such a way that 
\begin{equation}\label{ZUL1}
\Vert \mathcal K\Vert_{L^2((0,2)\times (0,T))} \le C,
\end{equation}
provided that  function $\psi$ belongs to some bounded set in $C^2[0,T].$ We give the proof of this fact in appendix.

Consider  function $w(t,x)=\int_{0}^T\mathcal K(t,\tau)u(\tau,x)d\tau$ and $r(t,x)=\int_{0}^T\mathcal K(t,\tau)\partial_\nu u(\tau,x)d\tau.$  Observe that
the function $w$ belongs to the space  $L^2((0,2)\times \Omega).$ Indeed 
$$
\Vert w\Vert^2_{L^2((0,2)\times \Omega)}=\int_0^2\Vert \int_{0}^T\mathcal K(t,\tau)u(\tau,x)d\tau\Vert^2_{L^2(\Omega)}dt\le \int_0^2 \int_{0}^T\Vert \mathcal K(t,\tau) u(\tau,x)\Vert^2_{L^2(\Omega)}d\tau dt
$$
$$\le \int_0^2  \int_{0}^T\vert\mathcal  K(t,\tau)\vert^2 d\tau dt\int_{0}^T \Vert u(\tau,\cdot)\Vert^2_{L^2(\Omega)}d\tau .
$$
Similarly,
by (\ref{ZUL1}) we have 
\begin{equation}\label{ZUL2}
\Vert r \Vert_{L^2((0,2)\times \Gamma)}\le C\Vert \partial_\nu u\Vert_{L^2((0,T)\times \Gamma)}.
\end{equation}
Function $w$ satisfies the equation
\begin{equation}\label{Ssimas}
L(x,D) w=-\partial_{t}^2w+A(x,D')w=0 \quad \mbox{in}\quad  (0,2)\times  \Omega,
\end{equation}
\begin{equation}
\quad w=0\quad \mbox{on}\quad (0,2)\times \Omega,\quad  \partial_\nu w=r\quad \mbox{on}\quad (0,2)\times \Gamma.
\end{equation}

So, if \begin{equation}\label{lsinus}\Vert u\Vert_{H^{1,2}(Q)} \le M,\end{equation} since 
$$
\Vert w(0,\cdot)\Vert_{H^2(\Omega)}+\Vert w(2,\cdot)\Vert_{H^2(\Omega)}\le C\Vert u\Vert_{L^2(0,T;H^2(\Omega))}
$$ we have \begin{equation}\label{ZUL22X}\Vert w\Vert_{H^2((0,2)\times \Omega)}\le C M.\end{equation}

Next we claim that there exists $\hat\theta>0$ and  constant $C$ , both are independent of $w$,  such that
\begin{equation}\label{LOX}\Vert w(1,\cdot)\Vert_{L^2(\Omega)}+\Vert \partial_{t}w(1,\cdot)\Vert_{L^2(\Omega)}\le C M\Vert r\Vert^{\hat \theta}_{L^2((0,2)\times \Gamma)}. \end{equation}

Indeed, let $\psi\in C^2(\bar \Omega)$ be a function such that $\nabla \psi \ne 0$ on $\Omega,$ and $\partial_\nu\psi<0$ on $\partial\Omega\setminus\Gamma, \psi>0$ on $\Omega.$ (For existence of such a function  see e.g. \cite{Im}.) Let $\phi(t,x)=\psi(x)-K (t-1)^2 $ with sufficiently large $K$  such that
\begin{equation}\label{sosiska}
\inf_{x\in \Omega}\phi(1,x)>\mbox{max}\{\Vert \phi(0,\cdot)\Vert_{C^0(\bar\Omega)},\Vert \phi(2,\cdot)\Vert_{C^0(\bar\Omega)}\}.\end{equation} 
We set  $\varphi(t,x)=e^{\lambda\phi(t,x)}$. By direct computations one can check that there exists $\lambda _0$ such that for all $\lambda\ge \lambda_0$ the function $\varphi(t,x)$ is pseudo-convex with respect to 
the symbol $-\xi_0^2-\sum_{k,j=1}^n a_{kj}(x) \xi_k\xi_j.$ 

By (\ref{sosiska}), for all $\lambda$ sufficiently large,  one can chose numbers $\varphi_I$ and $\varphi_{II}$ such that
\begin{equation}\label{son}
\varphi_I<\mbox{inf}_{x\in \Omega} \varphi(1,x), \quad\mbox{and}\,\, \varphi_I>\varphi_{II}>\mbox{max}\{\Vert \varphi(0,\cdot)\Vert_{C^0(\bar\Omega)},\Vert \varphi(2,\cdot)\Vert_{C^0(\bar\Omega)}\}.\end{equation} 
We fix some $\lambda>\lambda_0$ such that (\ref{son}) holds true.

Since function $\varphi$ is pseudo-convex with respect to 
the symbol $-\xi_0^2-\sum_{k,j=1}^n a_{kj}(x) \xi_k\xi_j$ the following Carleman estimate for the second order elliptic equation holds true (see {e.g. \cite{Im}):
\begin{eqnarray}\label{Ssosiska1}
s^3\Vert we^{s\varphi}\Vert^2_{L^2((0,2)\times \Omega)}+ s\Vert \nabla we^{s\varphi}\Vert^2_{L^2((0,2)\times \Omega)}+ \Vert L_+(x,D,s)(we^{s\varphi})\Vert^2_{L^2((0,2)\times \Omega)}\nonumber\\ \le C (s^3\Vert we^{s\varphi}(0,\cdot)\Vert^2_{L^2(\Omega)}+s\Vert( \nabla we^{s\varphi})(0,\cdot)\Vert^2_{L^2(\Omega)}\nonumber\\+s^3\Vert we^{s\varphi}(2,\cdot)\Vert^2_{L^2(\Omega)}+s\Vert( \nabla we^{s\varphi})(2,\cdot)\Vert^2_{L^2(\Omega)}+s\Vert \partial_\nu we^{s\varphi}\Vert^2_{L^2((0,2)\times \Gamma)})\quad \forall s\ge s_0, \end{eqnarray}
where 
$$L_+(x,D,s)=\partial_t^2+\sum_{k,j=1}^n a_{kj}\partial^2_{x_kx_j}+ s^2((\partial_t \varphi)^2+\sum_{k,j=1}^n a_{kj}\partial_{x_k}\varphi\partial_{x_j}\varphi).
$$
The inequality (\ref{Ssosiska1}) implies 
\begin{eqnarray}\label{SSsosiska1}
\frac{1}{s}\Vert \partial_t^2 w+\sum_{k,j=1}^n a_{kj}\partial^2_{x_kx_j}w\Vert^2_{L^2((0,2)\times \Omega)}
\le C (s^3\Vert we^{s\varphi}(0,\cdot)\Vert^2_{L^2(\Omega)}+s\Vert( \nabla we^{s\varphi})(0,\cdot)\Vert^2_{L^2(\Omega)}\nonumber\\
+s^3\Vert we^{s\varphi}(2,\cdot)\Vert^2_{L^2(\Omega)}+s\Vert( \nabla we^{s\varphi})(2,\cdot)\Vert^2_{L^2(\Omega)}+s\Vert \partial_\nu we^{s\varphi}\Vert^2_{L^2((0,2)\times \Gamma)})\quad \forall s\ge s_0.\end{eqnarray}
Let $\rho(t)\in C^\infty_0 [\frac 18, \frac{15}{16}]$ be such that $\rho=1$ on $(\frac 12,\frac32).$
By (\ref{SSsosiska1}) and (\ref{Ssosiska1}) we have 
\begin{eqnarray}\label{SSSsosiska1}
\frac{1}{s}\Vert \partial_t^2 (\rho w)+\sum_{k,j=1}^n a_{kj}\partial^2_{x_kx_j}(\rho w)\Vert^2_{L^2((0,2)\times \Omega)}
\le C (s^3\Vert we^{s\varphi}(0,\cdot)\Vert^2_{L^2(\Omega)}+s\Vert( \nabla we^{s\varphi})(0,\cdot)\Vert^2_{L^2(\Omega)}\nonumber\\
+s^3\Vert we^{s\varphi}(2,\cdot)\Vert^2_{L^2(\Omega)}+s\Vert( \nabla we^{s\varphi})(2,\cdot)\Vert^2_{L^2(\Omega)}+s\Vert \partial_\nu we^{s\varphi}\Vert^2_{L^2((0,2)\times \Gamma)})\quad \forall s\ge s_0.\end{eqnarray}

 Function $\rho w$ satisfies the zero Dirichlet boundary conditions  on $\partial((\frac 18, \frac{15}{16})\times \Omega).$ Applying the $L^2$ a priori estimates for elliptic operators we have
\begin{eqnarray}\label{SSSSsosiska1}
\frac{1}{s}\Vert \sum_{k,j=0}\partial^2_{x_kx_j}(\rho w)\Vert^2_{L^2((0,2)\times \Omega)}
\le C (s^3\Vert we^{s\varphi}(0,\cdot)\Vert^2_{L^2(\Omega)}+s\Vert( \nabla we^{s\varphi})(0,\cdot)\Vert^2_{L^2(\Omega)}\nonumber\\
+s^3\Vert we^{s\varphi}(2,\cdot)\Vert^2_{L^2(\Omega)}+s\Vert( \nabla we^{s\varphi})(2,\cdot)\Vert^2_{L^2(\Omega)}+s\Vert \partial_\nu we^{s\varphi}\Vert^2_{L^2((0,2)\times \Gamma)})\quad \forall s\ge s_0.\end{eqnarray}
Finally, combining estimates (\ref{SSSsosiska1}) and (\ref{SSSSsosiska1}) we obtain
\begin{eqnarray}\label{sosiska1}
s^3\Vert we^{s\varphi}\Vert^2_{L^2((0,2)\times \Omega)}+ s\Vert \nabla we^{s\varphi}\Vert^2_{L^2((0,2)\times \Omega)}+\frac 1s \Vert \sum_{j,k=0}^n \partial_{x_j}\partial_{x_k}(we^{s\varphi})\Vert^2_{L^2((\frac 12,\frac32)\times \Omega)}\nonumber\\ \le C (s^3\Vert we^{s\varphi}(0,\cdot)\Vert^2_{L^2(\Omega)}+s\Vert( \nabla we^{s\varphi})(0,\cdot)\Vert^2_{L^2(\Omega)}\nonumber\\+s^3\Vert we^{s\varphi}(2,\cdot)\Vert^2_{L^2(\Omega)}+s\Vert( \nabla we^{s\varphi})(2,\cdot)\Vert^2_{L^2(\Omega)}+s\Vert \partial_\nu we^{s\varphi}\Vert^2_{L^2((0,2)\times \Gamma)})\quad \forall s\ge s_0. \end{eqnarray}
 
Setting $\varphi_M>\Vert \varphi\Vert_{C^0([0,2]\times \bar\Omega)}$ from (\ref{sosiska1}) and the trace theorem we have
\begin{eqnarray}\label{Ksinus}
\frac 1s(\Vert w(1,\cdot)\Vert^2_{H^1( \Omega)}+ \Vert \partial_{t} w(1,\cdot)\Vert^2_{L^2( \Omega)})e^{s\varphi_I}\nonumber\\\le C e^{s\varphi_{II}}(\Vert w(0,\cdot)\Vert^2_{L^2(\Omega)}+\Vert( \nabla w)(0,\cdot)\Vert^2_{L^2(\Omega)}+\Vert w(2,\cdot)\Vert^2_{L^2(\Omega)}\nonumber\\+\Vert( \nabla w)(2,\cdot)\Vert^2_{L^2(\Omega)})+Ce^{s\varphi_M}\Vert \partial_\nu w\Vert^2_{L^2((0,2)\times \Gamma)} \quad\forall s\ge s_0. \end{eqnarray}
Dividing both sides of inequality (\ref{Ksinus}) by $e^{s\varphi_I}/s$ we obtain
\begin{eqnarray}\label{Ssinus}
\Vert w(1,\cdot)\Vert^2_{H^1( \Omega)}+ \Vert \partial_{t} w(1,\cdot)\Vert^2_{L^2( \Omega)}\nonumber\\\le C e^{-s\alpha}(\Vert w(0,\cdot)\Vert^2_{L^2(\Omega)}+\Vert( \nabla w)(0,\cdot)\Vert^2_{L^2(\Omega)}\nonumber\\+\Vert w(2,\cdot)\Vert^2_{L^2(\Omega)}+\Vert( \nabla w)(2,\cdot)\Vert^2_{L^2(\Omega)})+Ce^{s\beta}\Vert \partial_\nu w\Vert^2_{L^2((0,2)\times \Gamma)}\quad \forall s\ge s_0. \end{eqnarray}
Here $\alpha=\varphi_I-\varphi_{II}, \beta=\varphi_M-\varphi_{I}.$ By (\ref{son}) these numbers are strictly positive.
From (\ref{Ssinus}), using (\ref{ZUL22X}) we have
\begin{equation}\label{XU}
\Vert w(1,\cdot)\Vert^2_{H^1( \Omega)}+ \Vert \partial_{t} w(1,\cdot)\Vert^2_{L^2( \Omega)}\le C e^{-s\alpha}M^2+Ce^{s\beta}\Vert \partial_\nu w\Vert^2_{L^2((0,2)\times \Gamma)}. \end{equation}
Let 
\begin{equation}M^2\Vert \partial_\nu w\Vert^2_{L^2((0,2)\times \Gamma)} < \mbox{min}\{ e^{-s_0 (\alpha+\beta)},1\}.\end{equation}
Setting  in (\ref{XU})
$s=-a\log (M^2\Vert \partial_\nu w\Vert^2_{L^2((0,2)\times \Gamma)}), a=\frac{1}{\alpha+\beta}$ we obtain
\begin{eqnarray}\label{XOLk}
\Vert w(1,\cdot)\Vert^2_{H^1( \Omega)}+ \Vert \partial_{t} w(1,\cdot)\Vert^2_{L^2( \Omega)}\nonumber\\\le C (M^2\Vert \partial_\nu w\Vert^2_{L^2((0,2)\times \Gamma)})^{a\alpha}M^2+C(M^2\Vert \partial_\nu w\Vert^2_{L^2((0,2)\times \Gamma)})^{-a\beta}\Vert \partial_\nu w\Vert^2_{L^2((0,2)\times \Gamma)}\nonumber\\ \le CM^{1(1-\theta)}\Vert \partial_\nu w\Vert^{2\theta}_{L^2((0,2)\times \Gamma)}, \end{eqnarray}
where $\theta=\alpha/(\alpha+\beta).$
Hence the proof of estimate  (\ref{LOX}) is complete.

Next we need the analog of the estimate (\ref{XOLk}) for $t$ from $(0,1].$
We claim   that there exists a positive $\theta_0$ such that
\begin{equation}\label{LOX1}\Vert w(t,\cdot)\Vert_{L^2( \Omega)}+\Vert \partial_{t}w(t,\cdot)\Vert_{L^2( \Omega)}\le C M\Vert r\Vert^{t\theta_0}_{L^2((0,2)\times \Gamma)}\quad \forall t\in (0,\frac 12]. \end{equation}
In order to prove the estimate  (\ref{LOX1}) we  again use technique bases  on the Carleman estimates. Observe that the function $\tilde\varphi(t)=e^{\lambda t}-1$ for all $\lambda$ sufficiently large is strictly pseudoconvex with
respect to the symbol $-\xi_0^2-\sum_{k,j=1}^n a_{kj}(x) \xi_k\xi_j.$
Then the following Carleman estimate holds true (see {e.g. \cite{Im}):

\begin{eqnarray}\label{AKsosiska1}
s^3\Vert we^{s\tilde\varphi}\Vert^2_{L^2((0,1)\times \Omega)}+ s\Vert \nabla we^{s\tilde\varphi}\Vert^2_{L^2((0,1)\times \Omega)}+ \Vert (L(x,D)+(s\partial_t\tilde \varphi)^2)(we^{s\tilde\varphi})\Vert^2_{L^2((0,1)\times \Omega)}\nonumber\\ \le C (s^3\Vert we^{s\tilde\varphi}(0,\cdot)\Vert^2_{L^2(\Omega)}+s\Vert( \nabla we^{s\varphi})(0,\cdot)\Vert^2_{L^2(\Omega)}\nonumber\\+s^3\Vert we^{s\tilde\varphi}(1,\cdot)\Vert^2_{L^2(\Omega)}+s\Vert( \nabla we^{s\tilde\varphi})(1,\cdot)\Vert^2_{L^2(\Omega)})\quad \forall s\ge s_1. \end{eqnarray}
Repeating the arguments (\ref{SSsosiska1})-(\ref{SSSSsosiska1}) we obtain from (\ref{AKsosiska1}):

\begin{eqnarray}\label{Ksosiska1}
s^3\Vert we^{s\tilde\varphi}\Vert^2_{L^2((0,1)\times \Omega)}+ s\Vert \nabla we^{s\tilde\varphi}\Vert^2_{L^2((0,1)\times \Omega)}+\frac 1s \Vert \sum_{k,j=1}^n \partial^2_{x_kx_j}(we^{s\tilde\varphi})\Vert^2_{L^2((0,\frac 12)\times \Omega)}\nonumber\\ \le C (s^3\Vert we^{s\tilde\varphi}(0,\cdot)\Vert^2_{L^2(\Omega)}+s\Vert( \nabla we^{s\tilde\varphi})(0,\cdot)\Vert^2_{L^2(\Omega)}\nonumber\\+s^3\Vert we^{s\tilde\varphi}(1,\cdot)\Vert^2_{L^2(\Omega)}+s\Vert( \nabla we^{s\tilde\varphi})(1,\cdot)\Vert^2_{L^2(\Omega)})\quad \forall s\ge s_2. \end{eqnarray}

From (\ref{Ksosiska1}) and the trace theorem  for all $t$ in $(0,\frac 12]$ we have
\begin{eqnarray}\label{KKsinus}
\frac 1s(\Vert w(t,\cdot)\Vert^2_{H^1( \Omega)}+ \Vert \partial_{t} w(t,\cdot)\Vert^2_{L^2( \Omega)})e^{s\tilde\varphi(t)}\le C (\Vert w(0,\cdot)\Vert^2_{L^2(\Omega)}\\+\Vert( \nabla w)(0,\cdot)\Vert^2_{L^2(\Omega)}+e^{s\tilde\varphi(1)}(\Vert w(1,\cdot)\Vert^2_{L^2(\Omega)}+\Vert( \nabla w)(1,\cdot)\Vert^2_{L^2(\Omega)})\forall s\ge s_2.\nonumber \end{eqnarray}
Dividing both sides of inequality (\ref{Ksinus}) by $e^{s\tilde\varphi(t)}/s$  for some positive $\tilde\alpha $ and $\tilde\beta$ we obtain
\begin{eqnarray}\label{KSsinus}
\Vert w(t,\cdot)\Vert^2_{H^1( \Omega)}+ \Vert \partial_{t} w(t,\cdot)\Vert^2_{L^2( \Omega)}\nonumber\\\le C e^{-s\tilde\alpha t}(\Vert w(0,\cdot)\Vert^2_{L^2(\Omega)}+\Vert( \nabla w)(0,\cdot)\Vert^2_{L^2(\Omega)})\nonumber\\+e^{s\tilde\beta t}\Vert w(1,\cdot)\Vert^2_{L^2(\Omega)}+\Vert( \nabla w)(1,\cdot)\Vert^2_{L^2(\Omega)})\quad \forall s\ge s_2\,\,\mbox{and}\,\,\forall t\in (0,\frac 12]. \end{eqnarray}
 
From (\ref{KSsinus}), using (\ref{ZUL22X}) we have
\begin{equation}\label{XXU}
\Vert w(t,\cdot)\Vert^2_{H^1( \Omega)}+ \Vert \partial_{t} w(t,\cdot)\Vert^2_{L^2( \Omega)}\le C e^{-s\tilde\alpha t}M^2+Ce^{s\tilde\beta t}\Vert \partial_\nu w\Vert^2_{L^2((0,2)\times \Gamma)}\quad \forall s\ge s_2\,\,\,\,\forall t\in (0,\frac 12]. \end{equation}

Let 
\begin{equation}\label{lopukh}M^2\Vert \partial_\nu w\Vert^2_{L^2((0,2)\times \Gamma)} < \mbox{min}\{ e^{-s_0 (\alpha+\beta)},e^{-s_2 (\tilde \alpha+\tilde \beta)},1\}.\end{equation}
Setting  in (\ref{XXU})
$s=-\tilde a\ln (M^2\Vert \partial_\nu w\Vert^2_{L^2((0,2)\times \Gamma)}), \tilde a=\frac{1}{\tilde \alpha+\tilde\beta }$ we obtain
\begin{eqnarray}
\Vert w(t,\cdot)\Vert^2_{H^1( \Omega)}+ \Vert \partial_{t} w(t,\cdot)\Vert^2_{L^2( \Omega)}\nonumber\\\le C (M^2\Vert \partial_\nu w\Vert^2_{L^2((0,2)\times \Gamma)})^{\tilde a\tilde\alpha t}M^2+C(M^2\Vert \partial_\nu w\Vert^2_{L^2((0,2)\times \Gamma)})^{-\tilde a\tilde\beta t}\Vert \partial_\nu w\Vert^2_{L^2((0,2)\times \Gamma)}\nonumber\\ \le CM^2\Vert \partial_\nu w\Vert^{2t\theta_0}_{L^2((0,2)\times \Gamma)}\quad  \forall t\in (0,\frac 12],\end{eqnarray}
where $\theta_0=\alpha/(\alpha+\beta).$
Hence the proof of estimate  (\ref{LOX1}) is complete.

Finally, using (\ref{ZUL2}) and (\ref{LOX1}), we estimate the trace of the function $w$ at time moment $t=0:$
\begin{eqnarray}\label{lopukh1}
\Vert w(0,\cdot )\Vert_{L^2(\Omega)}\le \Vert \int_0^\frac 12 \partial_{t} wdt \Vert_{L^2(\Omega)} +\Vert w(\frac 12,\cdot )\Vert_{L^2(\Omega)} \nonumber\\\le \int_0^\frac 12C M\Vert r\Vert^{t\theta_0}_{L^2((0,2)\times \Gamma)}dt+C M\Vert r\Vert^{\theta_0}_{L^2((0,2)\times \Gamma)}=\nonumber\\
CM (\int_0^\frac 12 e^{t\log(\Vert r\Vert^{\theta_0}_{L^2((0,2)\times \Gamma)})}dt+\Vert r\Vert^{\theta_0}_{L^2((0,2)\times \Gamma)})\nonumber\\
=CM\left (\frac{1}{\log\Vert r\Vert^{\theta_0}_{L^2((0,2)\times \Gamma)}}(e^{\frac 12\log(\Vert r\Vert^{\theta_0}_{L^2((0,2)\times \Gamma)})}-1)+\Vert r\Vert^{\theta_0}_{L^2((0,2)\times \Gamma)}\right )\nonumber\\
\le C/\vert\log\Vert \partial_\nu u\Vert_{L^2((0,T)\times \Gamma)}.\vert.
\end{eqnarray}
Here in order to get the last inequality we used (\ref{ZUL2}). If inequality (\ref{lopukh}) fails the estimate (\ref{lopukh1}) is trivial.

Next we set  $\tilde w(t,x)=\int_{0}^T\mathcal K(t,\tau)\partial_\tau u(\tau,x)d\tau$ and $\tilde r(t,x)=\int_{0}^T\mathcal K(t,\tau)\partial_\nu \partial_\tau u(\tau,x)d\tau.$
Repeating the previous arguments  we obtain the estimate similar to (\ref{lopukh1}):
\begin{eqnarray}\label{lopukh2}
\Vert \tilde w(0,\cdot )\Vert_{L^2(\Omega)}
\le C/\vert\log\Vert \partial_\nu\partial_t u\Vert_{L^2((0,T)\times \Gamma)}\vert.
\end{eqnarray}
Now we estimate $L^2$ norm of the function $u$.   The a priori for the Schr\"odinger equation imply
\begin{equation}\label{lopukh3}\Vert u\Vert_{C([0,T];L^2(\Omega))}\le C\Vert u (0,\cdot)\Vert_{L^2(\Omega)}.\end{equation}
Next  we estimate $u(0,\cdot).$  We use the following notations
$$
w_\psi(t,x)=\int_{0}^T\mathcal K(t,\tau)u(\tau,x)d\tau, \quad \tilde w_\psi(t,x)=\int_{0}^T\mathcal K(t,\tau)\partial_\tau u(\tau,x)d\tau,
$$
where $\mathcal K$ is solution to the controllability problem (\ref{HUKK}). Short computations imply
$$
w_{e^{i\kappa t}}(0,x)=\int_{0}^T e^{i\kappa \tau} u(\tau,x)d\tau=- \frac 1{i\kappa}\int_{0}^T (e^{i\kappa \tau}-e^{i\kappa T})\partial_\tau u(\tau,x)d\tau-
\frac 1{i\kappa} (e^{-i\kappa T}-e^{i\kappa T}) u(0,x).
$$
Then, taking $\kappa$ such that $(e^{-i\kappa T}-e^{i\kappa T})\ne 0$ we obtain
\begin{equation}\label{lopukh4}
 u(0,x)=\frac{i\kappa}{(e^{-i\kappa T}-e^{i\kappa T})}\left (w_{e^{i\kappa t}}(0,x) +\frac 1{i\kappa}(\tilde  w_{e^{i\kappa t}}(0,x)-\tilde w_{e^{i\kappa T}}(0,x))\right ).
\end{equation}

Estimating the right hand side of (\ref{lopukh4}) using (\ref{lopukh1}) and (\ref{lopukh2}) we obtain 
$$
\Vert u (0,\cdot)\Vert_{L^2(\Omega)}\le C/\vert\log\Vert \partial_\nu u\Vert_{L^2((0,T)\times \Gamma)}\vert+C/\vert\log\Vert \partial_\nu\partial_t u\Vert_{L^2((0,T)\times \Gamma)}\vert
$$
This estimate and (\ref{lopukh3}) implies (\ref{lopukh5}).
$\blacksquare$
\section{Proof of Theorem \ref{4}.}

{\bf Proof.}
The  Volterra   operator $K$ constructed in the proof of the theorem \ref{3} is invertible. Therefore 
\begin{eqnarray}
\Vert \partial_\nu z\Vert_{H^1(0,T;L^2(\Gamma))}\le C\Vert\partial_\nu u\Vert_{H^2(0,T;L^2(\Gamma))};\nonumber\\\quad \Vert z\Vert_{H^{1,2}((0,T)\times \Omega)}+\Vert \partial_tz\Vert_{H^{1,2}((0,T)\times \Omega)}\le C(\Vert u\Vert_{H^{1,2}((0,T)\times \Omega)}+\Vert \partial_t u\Vert_{H^{1,2}((0,T)\times \Omega)}+\Vert \partial^2_t u\Vert_{H^{1,2}((0,T)\times \Omega)})\nonumber
\end{eqnarray}
and function $z$ satisfies the equation 
\begin{equation}
P(x,D)z=0\quad\mbox{in}\,(0,T)\times \Omega,\quad z\vert_{(0,T)\times\partial\Omega}=0, \quad z(0,\cdot)=-if.\nonumber
\end{equation}
By theorem \ref{2} we have
\begin{equation}
\Vert z\Vert_{ L^2((0,T)\times \Omega)}
\le \frac{C}{\vert \log \Vert \partial_\nu z\Vert_{H^1(0,T;L^2( \Gamma))}\vert}
\le \frac{\tilde C}{\vert \log \Vert \partial_\nu u\Vert_{H^2(0,T;L^2( \Gamma))}
\vert}.
\end{equation} 
Then the standard a priori estimate for the  Schr\"odinger   equation imply  (\ref{cocos}).
$\blacksquare$
\section{Appendix}
Let $\Psi(t,\tau)\in C^2([0,3]\times [0,T])$ be a function such that $\Psi=\psi.$ Consider the controllability problem  for the Schr\"odinger equation with locally distributed control $v$: 
\begin{eqnarray}\label{1HUKK} MK=i\partial_{\tau}K-\partial^2_{t} K =M\Psi+v \,\,\mbox{in} \quad  t\in (0,3), \,\tau\in(0,T),\,\, \mbox{supp}\, v\subset (2,3)\times (0,T),\nonumber\\
\quad K(0,\tau)=0\,\,\mbox{on}\,\, (0,T),
\quad  K(\cdot,0)= \Psi (\cdot,0), \quad K(\cdot,T)=\Psi (\cdot,T).\end{eqnarray} 
Then function $\Psi-K$ would solve the controllability problem (\ref{HUKK}).

In order to solve problem (\ref{1HUKK}) it suffice  to prove an observability estimate
\begin{equation}\label{1XXP}\Vert p\Vert_{L^2(G)}\le C(\Vert p\Vert_{L^2(G_1)}+\Vert Mp\Vert_{L^2(G)}),\end{equation}
where $G=(0,3)\times (0,T), G_1=(2,3)\times (0,T).$
The estimate (\ref{1XXP}) imply (see e.g. \cite{Z})  the existence of solution to problem (\ref{1HUKK}) such that
$$
\Vert K\Vert_{L^2(G)}\le C(\Vert \Psi (\cdot,T)\Vert_{L^2(0,3)}+\Vert \Psi (\cdot,0)\Vert_{L^2(0,3)}+\Vert M\Psi\Vert_{L^2(G)}).
$$
In order to prove (\ref{1XXP})  we use the following observability estimate obtained in \cite {B-P}: let $b_1, b_0\in L^\infty(G)$ there exist a constant $C$ such that
\begin{equation}\label{XP}\Vert p\Vert_{L^2(0,T;H^1(0,3))}\le C (\Vert Mp +b_1\partial_t p+b_0p\Vert_{L^2(G)}+\Vert p\Vert_{L^2(0,T;H^1(2,3))}).  \end{equation}
Consider the function  $\tilde p(\tau,t)=\rho(t) (\int_t^3 \rho(\tilde t) p(\tau,\tilde t)d\tilde t-\int_0^3 \rho(\tilde t) p(\tau,\tilde t)d\tilde t),$  where $\rho(t)\in C^\infty[0,3],$$ \rho\vert_{[0,2]}=1,$$ \rho\vert_{[\frac 52,3]}=0.$ We derive equation for the function $\tilde p$. Observe that
$$
M (\rho p)=\rho Mp+2\partial_t\rho \partial_tp+\partial^2_t\rho p.
$$
Integrating this equation on the interval $(t,3)$ we have
$$
i\partial_\tau \int_t^3\rho p d\tilde t+\partial_t (\rho p)+f=M\int_t^3\rho p d\tilde t=M(\int_t^3\rho p d\tilde t-\int_0^3\rho p d\tilde t),
$$
$$
f=\int_t^3\rho Mpd\tilde t+2 \int_t^3\partial_t\rho \partial_tp d\tilde t+\int_t^3\partial^2_t\rho pd\tilde t.
$$
For the function $f$ we have the following estimate
\begin{equation}\label{gordon}
\Vert f\Vert_{L^2(G)}\le C(\Vert Mp\Vert_{L^2(G)}+\Vert p\Vert_{L^2(G_1)}).
\end{equation}
Finally 
$$
M\tilde p=-\rho f-\partial^2_t\rho (\int_t^3\rho p d\tilde t-\int_0^3\rho p d\tilde t)+2\partial_t\rho p.
$$

Applying the estimate (\ref{XP}) to the  above equation and using (\ref{gordon}) we obtain 
$$
\Vert\rho p\Vert_{L^2(G)}\le C(\Vert  \int_0^3 \rho pd\tilde t\Vert_{L^2(0,T)}+ \Vert M p\Vert_{L^2(G)}+\Vert K\Vert_{L^2(G_1)}).
$$
This estimate implies 
\begin{equation}\label{kapitan}
\Vert p(0,\cdot)\Vert_{L^2(0,3)}\le C
\Vert p\Vert_{L^2(G)}\le C(\Vert  \int_0^3 \rho pd\tilde t\Vert_{L^2(0,T)}+ \Vert M p\Vert_{L^2(G)}+\Vert p\Vert_{L^2(G_1)}).
\end{equation}
This estimate and uniqueness compactness argument imply (\ref{1XXP}).
Indeed suppose that (\ref{kapitan}) fails. Then there exist a sequence $\{p_n\}$ such that $\Vert p_n(0,\cdot)\Vert_{L^2(0,3)}=1$ and
\begin{equation}\label{kapitan1}
\Vert p_n(0,\cdot)\Vert_{L^2(0,3)}\ge n (\Vert M p_n\Vert_{L^2(G)}+\Vert p_n\Vert_{L^2(G_1)}).
\end{equation}
By (\ref{kapitan}) and Holgrem's uniqueness theorem  one can take a converging subsequence, which we denote again $p_n$ such that $p_n(0,\cdot)=\sum_{n=0}^\infty b_{n,k} e_k(t)\rightarrow 0$ weakly in $L^2(0,3)$ and $e_k(t)=\sin(\frac{k\pi t}{T}).$  In particularly
\begin{equation}\label{kiparis}b_{n,k} \rightarrow 0\quad \mbox{as}\quad n\rightarrow +\infty.\end{equation}

 By (\ref{kapitan}) we have
\begin{equation}\label{blevotina}\Vert p_n(0,\cdot)\Vert_{L^2(0,3)}\le C
\Vert p_n\Vert_{L^2(G)}\le C\Vert  \int_0^3 \rho p_nd t\Vert_{L^2(0,T)}+ o(1)\quad \mbox{as}\quad n\rightarrow \infty.\end{equation} 

Since the  function $p_n$ is given by formula $p_n(\tau,t)=\sum_{k=1}^\infty e^{ik^2\pi^2\tau/T^2}  b_{n,k} e_k(t).$ Using this formula we obtain
$$
\int_0^3 p_n(t,\tau)\rho(t) dt=\sum_{k=1}^\infty e^{ik^2\pi^2\tau/T^2}b_{n,k} \int_0^3 \rho(t) e_k(t)dt.
$$
On the other hand 
$$
\Vert  \int_0^3 \rho(t) p_n(\tau,t)d t\Vert^2_{L^2(0,T)}\le C \sum_{k=1}^\infty (b_{n,k} \int_0^3 \rho(t) e_k(t)dt)^2.
$$
Taking into account that $\int_0^3 \rho(t) e_k(t)dt=O(\frac {1}{k})$  and  using (\ref{kiparis}) we obtain that
$$
\Vert  \int_0^3 \rho(t) p_n(\tau,t)d t\Vert^2_{L^2(0,T)}\le C \sum_{k=1}^\infty (b_{n,k} \int_0^3 \rho(t) e_k(t)dt)^2\rightarrow 0\quad \mbox{as}\quad n\rightarrow \infty
$$ we obtain the contradiction to (\ref{blevotina}).
  $\blacksquare$

\end{document}